\documentclass{article}
\usepackage{amsfonts, mathrsfs, amsmath}
\newtheorem{def1}{Definition}[section]
\newtheorem{pro}{Proposition}[section]
\newtheorem{rem}{Remark}[section]

\begin{document}
\newtheorem{thm}{Theorem}[section]
\title{Another Sierpinski object in \textbf{BFTS}}
\author {Rana Noor\thanks{rananoor4@gmail.com} and Sheo Kumar Singh\thanks{sheomathbhu@gmail.com}\\ {\em Department of Mathematics}\\ {\em Banaras Hindu University},\\ {\em Varanasi-221 005, India}} 
\date{}
\maketitle 
\abstract
We identify a Sierpinski object in the category of fuzzy bitopological spaces, which is different from the two Sierpinski objects in this category, obtained earlier by Khastgir and Srivastava.
\section{Introduction}
\qquad Giuli and Salbany \cite{GS} studied the category $\textbf{BTOP}$ of bitopological spaces and identified two Sierpinski objects namely the `quad' and the `triad' in $\textbf{BTOP}$. Subsequently, Khastgir and Srivastava \cite{KS} studied the category $\textbf{BFTS}$ of fuzzy bitopological spaces and identified two Sierpinski objects `$I^{2}$' and `$2I$' in $\textbf{BFTS}$ and showed that they behave in the same way in $\textbf{BFTS}$ as do the `quad' and the `triad' in $\textbf{BTOP}$ respectively.

In the present work, we obtain another Sierpinski object in $\textbf{BFTS}$ which is different from `$I^{2}$' and `$2I$'.
\section{Preliminaries}
\qquad For the category theoretic notions used here, \cite{AHS} may be referred. All subcategories are assumed to be full and replete.\\

Let $L$ be a frame with $0$ and $1$ being its least and largest elements respectively. For a given set $X$, let $\bar{0}$ and $\bar{1}$ denote the constant maps from $X$ to $L$ with values $0$ and $1$ respectively.

Given a set $X$, $L^{X}$ will denote the family of all maps $\mu:X\rightarrow L$ (called \textit{$L$-sets} or \textit{$L$-fuzzy sets}). $L^{X}$ is also a frame under the frame structure induced by that on $L$.

We recall some definitions which are used in this paper.
\begin{def1}{\rm \cite{Roda}}
A family $\tau\subseteq L^{X}$ is called an {\rm $L$-topology} on a set $X$, and the pair $(X,\tau)$ an {\rm $L$-topological space}, if $\tau$ is closed under arbitrary suprema and finite infima. Furthermore, a map $f:(X,\tau)\rightarrow (Y,\delta)$ between two $L$-topological spaces is called {\rm continuous} if $f^{\leftarrow}(\nu)\in\tau$, for every $\nu\in\delta$.
\end{def1}

By $L\textbf{-TOP}$, we shall denote the category of all $L$-topological spaces and their continuous maps.

If we take $L=I(=[0,1])$, then an $L$-topological space is known as \textit{fuzzy topological space} (cf. \cite{Chang}).

Let $\textbf{FTS}$ denote the category of all fuzzy topological spaces and their continuous maps.
\begin{def1}{\rm \cite{KS}}
A {\rm fuzzy bitopological space} is a triple $(X,\tau_{1},\tau_{2})$, where $X$ is a set and $\tau_{1}$, $\tau_{2}$ are fuzzy topologies on $X$. Furthermore, if for every distinct pair $x,y\in X$, there exists $\mu\in \tau_{1}\cup\tau_{2}$ such that $\mu(x)\neq \mu(y)$, then $(X,\tau_{1},\tau_{2})$ is called $T_{0}$ $($or $2T_{0})$.
\end{def1}
\begin{def1}{\rm \cite{KS}}
A map $f:(X,\tau_{1},\tau_{2})\rightarrow (Y,\delta_{1},\delta_{2})$ between fuzzy bitopological spaces is called {\rm bicontinuous} $($resp. {\rm biopen}$)$ if $f^{\leftarrow}(\nu)\in \tau_{i}$, for every $\nu\in\delta_{i}$ $($resp. if $f^{\rightarrow}(\mu)\in\delta_{i}$, for every $\mu\in\tau_{i})$, for $i=1,2$.
\end{def1}

Let $\textbf{BFTS}$ denote the category of all fuzzy bitopological spaces and their bicontinuous maps and $\textbf{BFTS}_{0}$ denote the subcategory of $\textbf{BFTS}$ whose objects are $T_{0}$-fuzzy bitopological spaces.\\

The notions of \textit{subspace}, \textit{homeomorphism} and \textit{embedding}, for a fuzzy bitopological space, are on expected lines.\\

The category $\textbf{BFTS}$ has initial structures (cf. \cite{Alka}).

\begin{rem}
Let $\mathscr F=\{f_{i}:X\rightarrow (Y_{i},\delta_{i},\delta_{i}')$ $|$ $i\in I\}$ be a family of maps, where $X$ is a set and $\{(Y_{i},\delta_{i},\delta_{i}')$ $|$ $i\in I\}$ is a family of fuzzy bitopological spaces. Then the fuzzy bitopology $(\Delta, \Delta')$ on $X$, which is initial with respect to the family $\mathscr F$, is the one for which $\Delta$ $($resp. $\Delta')$ is the fuzzy topology on $X$ having the subbase $\{f_{i}^{\leftarrow}(\mu)$ $|$ $\mu\in\delta_{i}, i\in I\}$ $($resp. $\{f_{i}^{\leftarrow}(\mu')$ $|$ $\mu'\in\delta_{i}', i\in I\})$.
\end{rem}
\begin{def1}{\rm \cite{KS}}
Given a family $\{(X_{i},\tau_{i}, \tau_{i}')$ $|$ $i\in I\}$ of fuzzy bitopological spaces, the initial fuzzy bitopology on $X$ $(=\displaystyle\prod_{i\in I}X_{i})$ with respect to the family of all projection maps $\{p_{i}:X\rightarrow (X_{i},\tau_{i}, \tau_{i}')$ $|$ $i\in I\}$ is called the {\rm product fuzzy bitopology}.
\end{def1}

Let $\mathscr C$ be a category and $\mathscr H$ some class of $\mathscr C$-morphisms.
\begin{def1}{\rm \cite{KS}}
A $\mathscr C$-object $X$ is called:
\begin{itemize}
\item $\mathscr H$-{\rm injective} if for every $e:Y\rightarrow Z$ in $\mathscr H$ and every $\mathscr C$-morphism $f:Y\rightarrow X$, there exists a $\mathscr C$-morphism $g:Z\rightarrow X$ such that $g\circ e=f$.
\item a {\rm cogenerator} in $\mathscr C$ if for every pair of distinct $f,g\in \mathscr C (Y,Z)$, there exists $h\in \mathscr C (Z,X)$ such that $h\circ f\neq h\circ g$.
\end{itemize}
\end{def1}
\textbf{Note}: In many familiar categories, $\mathscr H$ is usually taken to be the class of all embeddings in these categories, in which case the term `$\mathscr H$-injective' is shortened to just `injective'.
\begin{def1}{\rm \cite{KS}}
Let $\mathscr A$ be a class of $\mathscr C$-objects. We say that $\mathscr C$ is {\rm $\mathscr H$-cogenerated by $\mathscr A$} if every $\mathscr C$-object is an $\mathscr H$-subobject $(X$ is an $\mathscr H$-subobject of $Y$ if there is some $h:X\rightarrow Y$, with $h\in\mathscr H)$ of a product of objects in $\mathscr A$.
\end{def1}
\begin{def1}{\rm \cite{Mane}}
Given a category $\mathscr C$ of sets with structures, an object $S$ of $\mathscr C$ is called a {\rm Sierpinski object} if for every $X\in ob\mathscr C$, the family of all $\mathscr C$-morphisms from $X$ to $S$ is initial.
\end{def1}

The object $(L, \Delta)$ of $L\textbf{-Top}$, where $\Delta=\{\bar{0},id_{L},\bar{1}\}$, is easily verified to be $T_{0}$ and a Sierpinski object in $L\textbf{-Top}$ (as has been shown in \cite{Arun} for the case $L=[0,1]$).
\section{Sierpinski objects in BFTS}
\qquad The set $L=\{(a,b)\in I\times I$ $|$ $a+b\leq 1\}$ is a frame under the partial order $\leq$, defined as $(a,b)\leq (a',b')$ iff $a\leq a'$ and $b\geq b'$; the \textit{supremum} and the \textit{infimum} of $\{(a_{i},b_{i})\in L$ $|$ $i\in I\}$ in $(L,\leq)$ being $(\displaystyle\bigvee_{i\in I}a_{i}, \displaystyle\bigwedge_{i\in I}b_{i})$ and $(\displaystyle\bigwedge_{i\in I}a_{i}, \displaystyle\bigvee_{i\in I}b_{i})$ respectively and $(0,1)$ and $(1,0)$ being its least and largest elements.

\textbf{In this section, $L$ will denote this particular frame $\{(a,b)\in I\times I$ $|$ $a+b\leq 1\}$}.
\begin{rem}
Note that an $L$-set $\mu:X\rightarrow L$ can be identified with two maps $\mu_{1}:X\rightarrow [0,1]$ and $\mu_{2}:X\rightarrow [0,1]$ such that $\mu_{1}=p_{1}\circ\mu$ and $\mu_{2}=p_{2}\circ\mu$, where $p_{1}, p_{2}:L\rightarrow [0,1]$ are the two projection maps $($to the first and second `coordinates' respectively$)$. Thus an $L$-set and an $L$-topology are what some authors call as an {\rm intuitionistic fuzzy set (cf. \cite{Ata})} and an {\rm intuitionistic fuzzy topology (cf. \cite{Coker})}, respectively.\\
\end{rem}

Given any $L$-topological space $(X,\tau)$, it turns out that $\tau_{1}=\{p_{1}\circ\mu$ $|$ $\mu\in\tau\}$ and $\tau_{2}=\{\bar{1}-(p_{2}\circ\mu)$ $|$ $\mu\in\tau\}$ are fuzzy topologies on $X$ (cf. \cite{Chang}) and so $(X,\tau_{1},\tau_{2})\in ob\textbf{BFTS}$. Accordingly, for the $L$-Sierpinski space $(L,\Delta)$, also we get $(L,\Delta_{1},\Delta_{2})\in ob\textbf{BFTS}$, where $\Delta_{1}=\{\bar{0},p_{1},\bar{1}\}$ and $\Delta_{2}=\{\bar{0},\bar{1}-p_{2},\bar{1}\}$. 

We show in this section, that this object turns out to be a Sierpinski object in\textbf{ BFTS}.
\begin{rem}
We note that the two projection map $p_{1}, p_{2}:L\rightarrow [0,1]$ are such that $p_{1}\leq \bar{1}-p_{2}$. For if $(a,b)\in L$, then $a+b\leq 1$, whereby $a\leq 1-b$. Hence $p_{1}(a,b)\leq (\bar{1}-p_{2})(a,b)$.
\end{rem}

The following result is easy to verify.

\begin{pro}
For every $(X,\tau_{1},\tau_{2})\in ob \mathbf{BFTS}$ and for every $\mu\in \tau_{1}$ $($resp. $\mu\in \tau_{2})$, the map $h_{\mu}:(X,\tau_{1},\tau_{2})\rightarrow(L,\Delta_{1},\Delta_{2})$ defined as $h_{\mu}(x)=(\mu(x),0)$ $($resp. $h_{\mu}(x)=(0,1-\mu(x)))$ is a morphism in $\mathbf{BFTS}$, with $h_{\mu}^{\leftarrow}(p_{1})=\mu$ $($resp. $h_{\mu}^{\leftarrow}(\bar{1}-p_{2})=\mu)$.
\end{pro}
\begin{thm}
$(L,\Delta_{1},\Delta_{2})$ is a Sierpinski object in $\mathbf{BFTS}$.
\end{thm}
\textbf{Proof}: Let $(X,\tau_{1},\tau_{2})\in ob\textbf{BFTS}$ and $\mathscr{F}=\textbf{BFTS}((X,\tau_{1},\tau_{2}),(L,\Delta_{1},\Delta_{2}))$. Let $(Y,\delta_{1},\delta_{2})\in ob\textbf{BFTS}$ and $g:Y\rightarrow X$ be a map such that $f\circ g:(Y,\delta_{1},\delta_{2})\rightarrow (L,\Delta_{1},\Delta_{2})$ is bicontinuous for every $f\in\mathscr{F}$. We wish to show that $g$ is bicontinuous. Let $\mu\in \tau_{1}$. The bicontinuous map $h_{\mu}:(X,\tau_{1},\tau_{2})\rightarrow(L,\Delta_{1},\Delta_{2})$, described in Proposition $3.1$, is already in $\mathscr F$. Now, $g^{\leftarrow}(\mu)=g^{\leftarrow}(h_{\mu}^{\leftarrow}(p_{1}))=(h_{\mu}og)^{\leftarrow}(p_{1})$. Hence $g^{\leftarrow}(\mu)\in\delta_{1}$ (as $h_{\mu}\circ g$ is bicontinuous).

Similarly, for every $\mu\in \tau_{2}$, $g^{\leftarrow}(\mu)\in\delta_{2}$. So $g$ is bicontinuous. Thus $(L,\Delta_{1},\Delta_{2})$ is a Sierpinski object in $\textbf{BFTS}$. $\Box$
\begin{pro}
$(L,\Delta_{1},\Delta_{2})$ is $T_{0}$.\\
\end{pro}

We point out that, earlier, two interesting Sierpinski objects in $\textbf{BFTS}$ have been found by Khastgir and Srivastava in \cite{KS}, viz., (i) $(I^{2},\Pi_{1},\Pi_{2})$, where $I^{2}=I\times I$, $\Pi_{i}=\{\bar{0},\pi_{i},\bar{1}\}$, $i=1, 2$, and $\pi_{1},\pi_{2}:I^{2}\rightarrow I$ are the two projection maps and (ii) $(2I,\Omega_{1},\Omega_{2})$, where $2I=(I\times \{0\}) \cup (\{0\}\times I)$, $\Omega_{i}=\{\bar{0},q_{i},\bar{1}\}$, $i=1, 2$, and $q_{1},q_{2}:2I \rightarrow I$ are two maps defined as
\begin{equation*}
q_{1}(x)=
\begin{cases}
\alpha, & \text {if $x=(\alpha,0)\in I\times \{0\}$} \\
0, & \text{otherwise}
\end{cases}
\end{equation*} 
and 
\begin{equation*}
q_{2}(x)=
\begin{cases}
\alpha, & \text {if $x=(0,\alpha)\in \{0\}\times I$} \\
0, & \text{otherwise.}
\end{cases}
\end{equation*}
Furthermore, while both of these turned out to be cogenerators in $\textbf{BFTS}_{0}$, only $(I^{2},\Pi_{1},\Pi_{2})$ turned to be injective also. Thus it is natural to ask: in what respect(s) the `new found' Sierpinski object $(L,\Delta_{1},\Delta_{2})$ is similar to the Sierpinski objects of \cite{KS} in $\textbf{BFTS}$?\\

\begin{pro}
$(L,\Delta_{1},\Delta_{2})$ is a cogenerator in $\mathbf{BFTS}_{0}$.
\end{pro}
\textbf{Proof}: Consider any distinct pair $f,g:(X,\tau_{1},\tau_{2})\rightarrow (Y,\delta_{1},\delta_{2})$ of morphisms in $\textbf{BFTS}_{0}$. Then for some $x\in X$, $f(x)\neq g(x)$. As $(Y,\delta_{1},\delta_{2})$ is $T_{0}$, $\mu(f(x))\neq \mu(g(x))$ for some $\mu \in \delta_{1}\cup\delta_{2}$. If $\mu \in \delta_{1}$ (resp. $\mu \in \delta_{2}$), then by Proposition $3.1$, there exist a bicontinuous map $h_{\mu}:(Y,\delta_{1},\delta_{2})\rightarrow(L,\Delta_{1},\Delta_{2})$ defined as $h_{\mu}(y)=(\mu(y),0)$ (resp. $h_{\mu}(y)=(0,1-\mu(y))$). Clearly $h_{\mu}\circ f\neq h_{\mu}\circ g$. Thus $(L,\Delta_{1},\Delta_{2})$ is a cogenerator in $\textbf{BFTS}_{0}$. $\Box$
\begin{pro}
$(X,\tau_{1},\tau_{2})\in ob\mathbf{BFTS}_{0}$ iff $\mathscr F =\mathbf{BFTS}((X,\tau_{1},\tau_{2}), \\(L,\Delta_{1},\Delta_{2}))$ separates points of $X$.
\end{pro}
\textbf{Proof}: Let $(X,\tau_{1},\tau_{2})\in ob\textbf{BFTS}_{0}$ and $x,y\in X$ with $x\neq y$. Then $\mu(x)\neq\mu(y)$, for some $\mu\in\tau_{1}\cup\tau_{2}$. If $\mu\in\tau_{1}$ (resp. $\mu\in\tau_{2}$), then the bicontinuous map $h_{\mu}:(X,\tau_{1},\tau_{2})\rightarrow(L,\Delta_{1},\Delta_{2})$, described in Proposition $3.1$, is already in $\mathscr F$. Clearly $h_{\mu}(x)\neq h_{\mu}(y)$. Thus $\mathscr F$ separates points of $X$.

Conversely, let $\mathscr F$ separate points of $X$ and let $x,y\in X$ with $x\neq y$. Then $f(x)\neq f(y)$, for some $f\in \mathscr F$ and hence $f^{\leftarrow}(p_{1})\in\tau_{1}$ and $f^{\leftarrow}(\bar{1}-p_{2})\in\tau_{2}$. As $f(x)\neq f(y)$, either $p_{1}(f(x))\neq p_{1}(f(y))$ or $(\bar{1}-p_{2})(f(x))\neq (\bar{1}-p_{2})(f(y))$, showing that $(X,\tau_{1},\tau_{2})$ is $T_{0}$. $\Box$
\begin{pro}
$(L,\Delta_{1},\Delta_{2})$ $\mathscr H$-cogenerates $\mathbf{BFTS}_{0}$, where $\mathscr H$ is the class of all embeddings in $\mathbf{BFTS}_{0}$.
\end{pro}
\textbf{Proof}: Let $(X,\tau_{1},\tau_{2})\in ob \textbf{BFTS}_{0}$ and $\mathscr F= \textbf{BFTS}((X,\tau_{1},\tau_{2}), (L,\Delta_{1},\Delta_{2}))$. Define $e:(X,\tau_{1},\tau_{2})\rightarrow (L,\Delta_{1},\Delta_{2})^{\mathscr F}$ as $e(x)f=f(x)$, for every $x\in X$ and for every $f\in\mathscr F$. Let, for $f\in\mathscr F$, $\pi_{f}$ denote the $f$-th projection map. Then for every $x\in X$, $(\pi_{f}\circ e)(x)=\pi_{f}(e(x))=e(x)f=f(x)$ implying that $\pi_{f}\circ e=f$. Thus $e$ is bicontinuous. Let $x,y\in X$ with $x\neq y$. Then $\mu(x)\neq\mu(y)$, for some $\mu\in\tau_{1}\cup\tau_{2}$. If $\mu\in\tau_{1}$ (resp. $\mu\in\tau_{2}$), then the bicontinuous map $h_{\mu}:(X,\tau_{1},\tau_{2})\rightarrow(L,\Delta_{1},\Delta_{2})$, described in Proposition $3.1$, is already in $\mathscr F$. Clearly $h_{\mu}(x)\neq h_{\mu}(y)$, showing that $e(x)\neq e(y)$. Thus $e$ is injective. Let $\mu\in\tau_{1}$. Then for every $x\in X$, $e^{\rightarrow}(\mu)(e(x))=\vee\{\mu(x')$ $|$ $e(x')=e(x)\}=\mu(x)=p_{1}(\mu(x),0)=p_{1}(h_{\mu}(x))=p_{1}(e(x)h_{\mu})=p_{1}(\pi_{h_{\mu}}(e(x)))=(p_{1}\circ\pi_{h_{\mu}})(e(x))$, implying that $e^{\rightarrow}(\mu)=(p_{1}\circ\pi_{h_{\mu}})|_{e(X)}$. Similarly, if $\mu\in\tau_{2}$ then $e^{\rightarrow}(\mu)=((\bar{1}-p_{2})\circ\pi_{h_{\mu}})|_{e(X)}$. Thus $e:X\rightarrow e(X)$ is biopen, i.e., $e$ is an embedding. Hence $(L,\Delta_{1},\Delta_{2})$ $\mathscr H$-cogenerates $\mathbf{BFTS}_{0}$. $\Box$\\

For a fuzzy bitopological space $(X,\tau_{1},\tau_{2})$, $(X,\tau_{1}\vee\tau_{2})$ is a fuzzy topological space, where $\tau_{1}\vee\tau_{2}$ is the coarsest fuzzy topology on $X$ finer than $\tau_{1}$ and $\tau_{2}$.

Let $(X,\tau_{1},\tau_{2})$ be a fuzzy bitopological space. Put $pt(\tau_{1}\vee \tau_{2})=\{p:\tau_{1}\vee \tau_{2}\rightarrow I$ $|$ $p$ is a frame map$\}$. For $\mu\in \tau_{1}\vee \tau_{2}$, define $\mu^{s}:pt(\tau_{1}\vee \tau_{2})\rightarrow I$ as $\mu^{s}(p)=p(\mu)$, for every $p\in pt(\tau_{1}\vee \tau_{2})$. Then $\tau_{1}^{s}=\{\mu^{s}$ $|$ $\mu\in\tau_{1}\}$ and $\tau_{2}^{s}=\{\mu^{s}$ $|$ $\mu\in\tau_{2}\}$ are fuzzy topologies on $pt(\tau_{1}\vee \tau_{2})$ (cf. \cite{KS}).
\begin{def1}{\rm \cite{KS}}
A fuzzy bitopological space $(X,\tau_{1},\tau_{2})$ is called bisober if $\eta_{X}:(X,\tau_{1},\tau_{2})\rightarrow (pt(\tau_{1}\vee \tau_{2}),\tau_{1}^{s},\tau_{2}^{s})$, defined as $\eta_{X}(x)(\mu)=\mu(x)$, for every $x\in X$ and for every $\mu\in \tau_{1}\vee \tau_{2}$, is bijective.\\
\end{def1}

In \cite{KS}, both $(I^{2},\Pi_{1},\Pi_{2})$ and $(2I,\Omega_{1},\Omega_{2})$ are shown to be bisober.
\begin{pro}
$(L,\Delta_{1},\Delta_{2})$ is bisober.
\end{pro}
\textbf{Proof}: We show that $\eta_{L}: (L,\Delta_{1}, \Delta_{2})\rightarrow (pt(\Delta_{1}\vee \Delta_{2}), \Delta_{1}^{s},\Delta_{2}^{s})$ is bijective. The injectivity of $\eta_{L}$ easily follows from the fact that $(L,\Delta_{1},\Delta_{2})$ is $T_{0}$. Now we show that $\eta_{L}$ is surjective. Let $p\in pt(\Delta_{1}\vee \Delta_{2})$. Then $p:\Delta_{1}\vee \Delta_{2}\rightarrow I$, being a frame map, is order preserving. So if $p(p_{1})=\alpha$ and $p(\bar{1}-p_{2})=\beta$, then $\alpha\leq\beta$. Hence $\alpha+1-\beta \leq 1$, implying that $(\alpha, 1-\beta)\in L$. Clearly $\eta_{L}(\alpha, 1-\beta)=p$. Thus $\eta_{L}$ is surjective and hence $(L,\Delta_{1},\Delta_{2})$ is bisober. $\Box$
\begin{pro}
$(L,\Delta_{1},\Delta_{2})$ is not injective in $\mathbf{BFTS}_{0}$.
\end{pro}
\textbf{Proof}: Consider the identity map $id:(L,\Delta_{1},\Delta_{2})\rightarrow (L,\Delta_{1},\Delta_{2})$, which is clearly an extremal monomorphism. Define $e:(L,\Delta_{1},\Delta_{2})\rightarrow (I^{2},\Pi_{1},\Pi_{2})$ as $e(a,b)=(a,1-b)$. It is easy to see that $e$ is bicontinuous and injective. Also, for $(a,b)\in L$, $(e^\rightarrow (p_{1}))e(a,b)=\bigvee \{p_{1}(x,y)$ $|$ $e(x,y)=e(a,b)\}=p_{1}(a,b)=a$ and $\pi_{1}(e(a,b))=\pi_{1}(a,1-b)=a$. Thus, $e^\rightarrow (p_{1})=\pi_{1}|_{e(L)}$, whereby $e^\rightarrow (p_{1})\in \Pi_{1}|_{e(L)}$. Similarly, as $e^\rightarrow (\bar{1} -p_{2})=\pi_{2}|_{e(L)}$, $e^\rightarrow (\bar{1} -p_{2})\in \Pi_{2}|_{e(L)}$. Thus, $e: L\rightarrow e(L)$ is biopen. Hence, $e$ is an embedding. We show that $e$ is also an epimorphism in $\textbf{BFTS}_{0}$. Consider any distinct pair $f,g:(I^{2},\Pi_{1},\Pi_{2})\rightarrow (Y,\delta_{1},\delta_{2})$ of morphisms in $\textbf{BFTS}_{0}$. Then for some $x\in I^{2}$, $f(x)\neq g(x)$. As $Y$ is $T_{0}$, there exists $\mu\in \delta_{1} \cup\delta_{2}$ such that $\mu(f(x))\neq \mu (g(x))$, i.e. $f^{\leftarrow}(\mu)\neq g^{\leftarrow}(\mu)$. If $\mu\in \delta_{1}$ (resp. $\mu\in \delta_{2}$), then $f^{\leftarrow}(\mu), g^{\leftarrow}(\mu) \in \Pi_{1}$ (resp. $f^{\leftarrow}(\mu), g^{\leftarrow}(\mu) \in \Pi_{2}$). Now $(1/2,1/2)\in L$ and $f^{\leftarrow}(\mu)(1/2,1/2)\neq g^{\leftarrow}(\mu)(1/2,1/2)$. This implies that $(f\circ e)(1/2,1/2)\neq (g\circ e)(1/2,1/2)$, whereby $f\circ e\neq g\circ e$. Thus $e$ is an epimorphism.

Now if there exists a morphism $h:(I^{2},\Pi_{1},\Pi_{2})\rightarrow (L,\Delta_{1},\Delta_{2})$ in $\textbf{BFTS}_{0}$ such that $h\circ e=id$, then, as $id$ is an extremal monomorphism, $e$ will have to be an isomorphism, which clearly is not possible. Thus $(L,\Delta_{1},\Delta_{2})$ cannot be injective in $\textbf{BFTS}_{0}$. $\Box$\\

\begin{rem}
The above result shows that $(L,\Delta_{1},\Delta_{2})$ and $(I^{2},\Pi_{1},\Pi_{2})$ are different.
\end{rem}

In our last result, we shall use the following easy-to-verify result.
\begin{pro}
\begin{enumerate}
\item If $f:(X,\tau_{1},\tau_{2})\rightarrow (Y,\delta_{1},\delta_{2})$ is a homeomorphism in $\mathbf{BFTS}$ then $f:(X,\tau_{1}\vee\tau_{2})\rightarrow (Y,\delta_{1}\vee\delta_{2})$ is a homeomorphism in $\mathbf{FTS}$.
\item If $f:(X,\tau)\rightarrow (Y,\delta)$ is a homeomorphism in $\mathbf{FTS}$ then $f^{\leftarrow}:\delta\rightarrow\tau$ is bijective.
\end{enumerate}
\end{pro}
\begin{thm}
The fuzzy bitopological spaces $(2I,\Omega_{1},\Omega_{2})$ and $(L,\Delta_{1},\Delta_{2})$ are not homeomorphic.
\end{thm}
\textbf{Proof}: Consider the fuzzy topological spaces $(2I,\Omega_{1}\vee\Omega_{2})$ and $(L,\Delta_{1}\vee\Delta_{2})$. It is clear that $q_{1}\wedge q_{2}=\bar{0}$, so $\Omega_{1}\vee\Omega_{2}=\{\bar{0}, q_{1}, q_{2}, q_{1}\vee q_{2},\bar{1}\}$. As $p_{1}\leq (\bar{1}-p_{2})$, $\Delta_{1}\vee\Delta_{2}=\{\bar{0}, p_{1}, \bar{1}-p_{2}, \bar{1}\}$. This shows that the number of elements in $\Omega_{1}\vee\Omega_{2}$ and $\Delta_{1}\vee\Delta_{2}$ are not same. Hence there cannot exist any bijection between $\Omega_{1}\vee\Omega_{2}$ and $\Delta_{1}\vee\Delta_{2}$. Thus $(2I,\Omega_{1}\vee\Omega_{2})$ and $(L,\Delta_{1}\vee\Delta_{2})$ are not homeomorphic and hence $(2I,\Omega_{1},\Omega_{2})$ and $(L,\Delta_{1},\Delta_{2})$ are also not homeomorphic. $\Box$\\\\
\textbf{Acknowledgement}: The authors would like to thank Prof. A.K. Srivastava for providing counsel in the preparation of this paper. The authors (RN and SKS) would also like to respectively thank the {\it University Grants Commission} (New Delhi, India) and the {\it Council of Scientific \& Industrial Research} (New Delhi, India) for financial supports through their respective Senior Research Fellowships.

\end{document}